\documentclass[12pt]{amsart}
\usepackage{enumerate}

\begin{document}
\numberwithin{equation}{section}

\def\1#1{\overline{#1}}
\def\2#1{\widetilde{#1}}
\def\3#1{\widehat{#1}}
\def\4#1{\mathbb{#1}}
\def\5#1{\frak{#1}}
\def\6#1{{\mathcal{#1}}}

\def\C{{\4C}}
\def\R{{\4R}}
\def\N{{\4N}}
\def\Z{{\4Z}}

\title[Totally geodesic discs in bounded symmetric domains]
{Totally geodesic discs in bounded symmetric domains}
\author[S.-Y. Kim \& A. Seo]{Sung-Yeon Kim and Aeryeong Seo}
\address{S.-Y. Kim: Center for Complex Geometry, Institute for Basic Science, 55 Expo-ro, Yuseong-gu,
Daejeon 34126, South Korea }
\email{sykim8787@ibs.re.kr}

\address{A. Seo: Department of Mathematics,
Kyungpook National University,
Daegu 41566, Republic of Korea}%
\email{aeryeong.seo@knu.ac.kr}

%\thanks{*This research was supported by Basic Science Research Program through the National Research Foundation of Korea(NRF) funded by the Ministry of  Science, ICT and Future Planning(grant number NRF-2015R1A2A2A11001367)}
\maketitle
%\tableofcontents

%\def\Label#1{\label{#1}{\bf (#1)}~}
\def\Label#1{\label{#1}}

% Standard sets

\def\cn{{\C^n}}
\def\cnn{{\C^{n'}}}
\def\ocn{\2{\C^n}}
\def\ocnn{\2{\C^{n'}}}

% Abbreviations

\def\dist{{\rm dist}}
\def\const{{\rm const}}
\def\rk{{\rm rank\,}}
\def\id{{\sf id}}
\def\aut{{\sf aut}}
\def\Aut{{\text{Aut}}}
\def\CR{{\rm CR}}
\def\GL{{\sf GL}}
\def\Re{{\sf Re}\,}
\def\Im{{\sf Im}\,}
\def\span{\text{\rm span}}
\def\mult{\text{\rm mult\,}}
\def\reg{\text{\rm reg\,}}
\def\ord{\text{\rm ord\,}}
\def\hot{\text{\rm HOT\,}}

\def\codim{{\rm codim}}
\def\crd{\dim_{{\rm CR}}}
\def\crc{{\rm codim_{CR}}}

\def\eps{\varepsilon}
\def\d{\partial}
\def\a{\alpha}
\def\b{\beta}
\def\g{\gamma}
\def\G{\Gamma}
\def\D{\Delta}
\def\Om{\Omega}
\def\k{\kappa}
\def\l{\lambda}
\def\L{\Lambda}
\def\z{{\bar z}}
\def\w{{\bar w}}
\def\Z{{\1Z}}
\def\t{\tau}
\def\th{\theta}

\emergencystretch15pt
\frenchspacing

\newtheorem{Thm}{Theorem}[section]
\newtheorem{Cor}[Thm]{Corollary}
\newtheorem{Pro}[Thm]{Proposition}
\newtheorem{Lem}[Thm]{Lemma}
\newtheorem{Prob}[Thm]{Problem}

\theoremstyle{definition}\newtheorem{Def}[Thm]{Definition}

\theoremstyle{remark}
\newtheorem{Rem}[Thm]{Remark}
\newtheorem{Exa}[Thm]{Example}
\newtheorem{Exs}[Thm]{Examples}

\def\bl{\begin{Lem}}
\def\el{\end{Lem}}
\def\bp{\begin{Pro}}
\def\ep{\end{Pro}}
\def\bt{\begin{Thm}}
\def\et{\end{Thm}}
\def\bc{\begin{Cor}}
\def\ec{\end{Cor}}
\def\bd{\begin{Def}}
\def\ed{\end{Def}}
\def\br{\begin{Rem}}
\def\er{\end{Rem}}
\def\be{\begin{Exa}}
\def\ee{\end{Exa}}
\def\bpf{\begin{proof}}
\def\epf{\end{proof}}
\def\ben{\begin{enumerate}}
\def\een{\end{enumerate}}
\def\beq{\begin{equation}}
\def\eeq{\end{equation}}

\begin{abstract}
In this paper, we characterize $C^2$-smooth totally geodesic isometric embeddings $f\colon \Omega\to\Omega'$ between bounded symmetric domains $\Omega$ and $\Omega'$ which extend $C^1$-smoothly over some open subset in the Shilov boundaries and 
have nontrivial normal derivatives on it.
In particular, if $\Omega$ is irreducible, there exist totally geodesic bounded symmetric subdomains $\Omega_1$ and $\Omega_2$ of $\Omega'$ such that $f = (f_1, f_2)$ maps into $\Omega_1\times \Omega_2\subset \Omega$ where $f_1$ is holomorphic and $f_2$ is anti-holomorphic totally geodesic isometric embeddings. 
If $\text{rank}(\Omega')<2\text{rank}(\Omega)$, then either $f$ or $\bar f$ is a standard holomorphic embedding.
\end{abstract}

\section{Introduction}
Totally geodesic maps between Riemannian manifolds are one of the most important examples among harmonic maps. In \cite{Siu80, Siu81}, Siu showed that any harmonic map from a compact K\"ahler manifold is either holomorphic or anti-holomorphic if the rank of the differential of the map is at least $4$ at some point and if the curvature tensor of the target manifold is strongly negative for his celebrated strong rigidity theorem. Moreover, applying his method to wider classes of complex manifolds, he could show that the same result holds for bounded symmetric domains. 

On the other hand if complex manifolds are noncompact, the behavior of harmonic maps becomes more complicated as shown by Li-Ni in \cite{LN00}: even for the unit balls one needs extra conditions for harmonic maps to be holomorphic or anti-holomorphic. For research in this direction, we refer the reader to \cite{LS07} for strongly pseudoconvex domains and \cite{Xiao21} for bounded symmetric domains of classical type.

In this paper we examine complex-analyticity of totally geodesic isometric embeddings between bounded symmetric domains with respect to their K\"ahler-Einstein metrics.
Greene-Krantz (\cite{GK82}) proved that any isometry between strongly pseudoconvex domains with respect to their Bergman or K\"ahler-Einstein metrics is either holomorphic or anti-holomorphic and we refer the reader to \cite{Wu88} for the corresponding result for irreducible bounded symmetric domains. 

Let us now fix some notations.  For a bounded symmetric domain $\Omega$, denote by $g_\Omega$ the canonical K\"ahler-Einstein metric on $\Omega$ normalized so that minimal discs have constant Gaussian curvature $-2$. Note that $g_\Omega$ is a constant multiple of its Bergman metric.
We will denote by $\Aut(\Omega)$ the set of holomorphic diffeomorphisms of $\Omega$ onto itself and let
$\Delta=\{z\in \mathbb C : |z|<1\}$ be the unit disc in $\mathbb C$.

%\bd
%Let $\Omega$, $\Omega'$ be bounded symmetric domains and let $f\colon\Omega\to\Omega'$ be a $C^1$ map. $f$ is called an isometry with respect to the Bergman metric if there exists a constant $\lambda>0$ such that $f^*(d_{\Omega'}^2)=\lambda d_\Omega^2$. If $f$ is $C^2$-smooth, then $f$ is said to be totally geodesic if the second fundamental form of $f$ is trivial.
%\ed
In this paper we prove the following results: 
\bt\label{main thm}
Let $\Omega$ be a bounded symmetric domain and let $f\colon(\Delta, \lambda g_\Delta)\to (\Omega,g_\Omega)$ be a $C^2$-smooth totally geodesic isometric embedding for some $\lambda>0$.
Suppose that $f$ extends $C^1$ up to an open set $U\subset \partial\Delta$. Suppose further that the radial derivative of $f$ is nontrivial on $U$. Then there exists a totally geodesic polydisc $\Delta^r\subset \Omega$ such that $f(\Delta)\subset \Delta^r$ and if we express $f$ as $(f_1,\ldots,f_r)\colon\Delta\to \Delta^r$, then  each $f_i$ belongs to $\Aut(\Delta)$ or $\overline{\Aut(\Delta)}$.
\et

For Hermitian symmetric spaces of compact type $X$ and $X'$, a holomorphic map $F\colon X\to X'$ is called a {\it standard embedding} if there exists a characteristic subspace $X''\subset X'$ with $\text{rank}(X'') = \text{rank}(X)$  such that $F(X)\subset X''$ and $F\colon X\to X''$ is a totally geodesic isometric embedding with respect to their canonical K\"ahler-Einstein metric. For bounded symmetric domains $\Omega$ and $\Omega'$, we say that a holomorphic map $f\colon\Omega\to \Omega'$ is a {\it standard embedding} if $f$ is a restriction of a standard embedding $F\colon X\to X'$ to $\Omega$ where $X$ and $X'$ are the compact duals of $\Omega$ and $\Omega'$ respectively.

\bt\label{main cor}
Let $(\Omega,g_\Omega),~ (\Omega',g_{\Omega'})$ be bounded symmetric domains with normalized K\"ahler-Einstein metric and let $f\colon (\Omega, \lambda g_\Omega)\to (\Omega', g_{\Omega'})$ be a $C^2$-smooth totally geodesic isometric embedding for some $\lambda>0$. Suppose that $f$ extends $C^1$ up to an open set $U$ of the Shilov boundary of $\Omega$ and the radial derivative of $f$ is nontrivial on $U$. Suppose further that $\Omega$ is irreducible. Then there exist totally geodesic bounded symmetric subdomains $\Omega_1,\Omega_2$ of $\Omega'$ such that $f=(f_1,f_2)\colon\Omega\to \Omega_1\times\Omega_2$, where, unless they are constant, $f_1$ is holomorphic and $f_2$ is anti-holomorphic totally geodesic isometric embedding. In particular, if
$${\rm rank}(\Omega')<2~{\rm rank}(\Omega),$$ 
then either $f$ or $\bar f$ is a standard holomorphic embedding . 
\et

{\bf Acknowledgement} The first author was supported by the Institute for Basic Science (IBS-R032-D1-2021-a00). The second author was supported by Basic Science Research Program through the National Research Foundation
of Korea (NRF) funded by the Ministry of Education (NRF-2019R1F1A1060175).

\section{Preliminaries}
In this section we collect background information and known results. See \cite{Wolf72, Mok86} for more details.

A bounded domain $\Omega$ is called symmetric if for each $p\in \Omega$, there exists a holomorphic automorphism $I_p$ such that $I_p^2$ is the identity map of $\Omega$ which has $p$ as an isolated fixed point. 
Bounded symmetric domains are homogeneous complex manifolds and their Bergman metrics are K\"ahler-Einstein with negative holomorphic sectional curvatures.
All Hermitian symmetric spaces of non-compact type can be realized as convex bounded symmetric domains by the Harish-Chandra realizations.
Any Hermitian symmetric space of noncompact type can be canonically embedded into a Hermitian symmetric space of compact type, called the compact dual, by Borel embedding.

Throughout this article $M(p,q; \mathbb C)$ denotes the set of $p\times q$ matrices with complex coefficients.
The set of irreducible Hermitian symmetric spaces of non-compact type consists of four classical types and two exceptional types. We list the irreducible bounded symmetric domains of classical type which are the Harish-Chandra realizations of them as follows;
$$
\Omega^I_{p,q} := \left\{ Z \in M(p,q;\mathbb C) : I - \overline Z^t Z>0 \right\}, \quad 1\leq p\leq q;
$$
$$\Omega^{II}_n := \left\{ Z\in \Omega^I_{n,n} : Z^t = -Z \right\}, \quad n\geq 2;
$$
$$
\Omega^{III}_n := \left\{ Z\in \Omega^I_{n,n} : Z^t = Z \right\}, \quad n\geq 1;
$$
$$
\Omega^{IV}_n:= \left\{ (z_1, \ldots, z_n) \in \mathbb C^n : || z||^2 < 2, || z||^2 <1 + \bigg| \frac{1}{2} \sum_{j=1}^n z_j^2 \bigg|^2 \right\}, \quad n\geq 3;
$$
There are two exceptional type $\Omega^V_{16}$ and $\Omega^{VI}_{27}$.

Let $S_{p,q}^{I}$, $S_n^{II} $, $S_n^{III}$, $S_n^{IV}$, $S^{V}$ and $S^{VI}$ be generic norms of the corresponding domains which are given by
$$S_{p,q}^{I}(Z,\overline Z) = \det(I_r -ZZ^*) 
\quad \text{ for }  Z\in M(p,q;\mathbb C), $$
$$S_n^{II} (Z,\overline Z) = s_n^{II}(Z) 
\quad\text{ for } Z\in \{Z\in M(n,n;\mathbb C):Z^t=-Z\},$$
$$S_n^{III} (Z,\overline Z)= \det(I_n -ZZ^*) 
\quad\text{ for } Z\in \{Z\in M(n,n;\mathbb C):Z^t=Z\},$$
$$S_n^{IV}(Z,\overline Z) = 1-2 ZZ^* + \left| ZZ^t \right|^2 
\quad\text{ for } Z\in {\mathbb C}^n,$$ 
for the classical domains where $\det(I_n -ZZ^*) = s_n^{II}(Z)^2$ for some polynomial 
$s_n^{II}(Z)$ and $Z\in \{Z\in M(n,n;\mathbb C):Z^t=-Z\}$.
The Bergman kernel $K(z,z)$ of $\Omega$ can be
expressed by 
$$
c_1N_\Omega(z)^{-c_2}
$$ 
for some constant $c_1, c_2>0$ where $N_\Omega$ is the generic norm of $\Omega$. 
Therefore the Bergman metric is given by $$-c_2\sum \frac{\partial^2}{\partial z_j\partial \bar z_k} \log N_\Omega(z) dz_j\otimes d\bar z_k.$$
It is well known that its Ricci curvature is a negative constant.
\bt[Polydisc Theorem]\label{polydisc theorem}
Let $\Omega$ be a bounded symmetric domain and $g_\Omega$ be its Bergman metric. Let $X_\Omega$ be the compact dual of $\Omega$ and $g_c$ be its K\"ahler-Einstein metric.
There exists a totally geodesic complex submanifold $D$ of $(\Omega, g_\Omega)$ such that 
$(D, g_\Omega\big|_D)$ is holomorphically isometric to a Poincar\'e polydisc $(\Delta^r, \rho)$ and $$
\Omega = \bigcup_{\gamma\in K} \gamma D
$$
where $K$ denotes an isotropy subgroup of $\text{Aut}(\Omega)$.
Moreover, there exists a totally geodesic complex submanifold $S$ of $(X,g_c)$
containing $D$ as an open subset such that $(S, g_c\big|_S)$ is isometric to a polysphere 
$((\mathbb P^1)^r, \rho_c)$ equipped with a product Fubini-Study metric $\rho_c$.
\et

The dimension of $D$ in Theorem \ref{polydisc theorem} is called the {\it rank} of $\Omega$ and it is given by Table \ref{rank} for each irreducible bounded symmetric domain.
\begin{table}[ht]
\caption{rank of $\Omega$}\label{rank}
\begin{tabular}{c|c|c|c|c|c|c}
$\Omega$& $\Omega_{p,q}^I (p\leq q)$& $\Omega^{II}_n$& $\Omega^{III}_n$& $\Omega^{IV}_n$ & $\Omega^V_{16}$&$\Omega_{27}^{VI}$  \\[4pt]\hline
&&&&&&\\[-10pt]
rank & $p$ & $\left[ \frac{n}{2}\right]$& $n$ & $2$ & $2$ & $3$
\end{tabular}
\end{table}

A totally geodesic disc $\Delta$ which can be expressed by $\{ (z,0,\ldots,0): |z|<1\}\subset \Delta^r \cong D$ is called a {\it minimal disc}. We say that the vector $v\in T\Omega$ is of {\it rank} $k$ if the minimal polydisc tangential to $v$ is $k$-dimensional.
For a nonzero rank one vector $v\in T_p\Omega$, define
$$\mathcal N_{[v]}:=\{ w\in T_p\Omega: R(v,\bar v, w, \bar w)=0\},$$
where $R$ is the curvature tensor of $g_\Omega$.
Remark that any vector $v\in T\Omega$ can be expressed as a linear combination of rank one vectors by the Polydisc Theorem.
For a nonzero vector $v=\sum_j c_jv_j$ with rank one vectors $v_j$, define
$$\mathcal N_{[v]}:=\bigcap_j\mathcal N_{[v_j]}.$$
For a totally geodesic polydisc $\Delta^k\subset \Omega$, denote by $(\Delta^k)^\perp$ the totally geodesic subdomain of $\Omega$, called a characteristic subdomain of $\Omega$, such that 
$$T_p (\Delta^k)^\perp=\bigcap_{[v]\in \mathbb PT_p\Delta^k} \mathcal N_{[v]}.$$ 
For totally geodesic polysphere $S\cong (\mathbb P^1)^k$, we denote by $S^\perp$ the compact dual of $(\Delta^k)^\perp$ canonically embedded into the compact dual of $\Omega$ containing $(\Delta^k)^\perp$.
For each irreducible bounded symmetric domain $\Omega$, $(\Delta^k)^\perp$ is given by Table \ref{characteristic subdomains} below for each $k$ with $1\leq k\leq \text{rank}(\Omega)$ and the canonical embedding into $\Omega$.

\begin{table}[ht]\caption{Characteristic subdomains}\label{characteristic subdomains}
\begin{tabular}{c|c|c|c|c|c|c}
$\Omega$& $\Omega_{p,q}^I (p\leq q)$& $\Omega^{II}_n$& $\Omega^{III}_n$& $\Omega^{IV}_n$ & $\Omega^V_{16}$&$\Omega_{27}^{VI}$  \\[4pt]\hline 
&& &&&&\\[-8pt]
$(\Delta^k)^\perp$ & $\Omega^I_{p-k, q-k}$ & $\Omega^{II}_{n-2k}$& $\Omega_{n-k}^{III}$ & $\Delta(k=1)$ & $\mathbb B^5(k=1)$ & $\Omega_8^{IV} (k=2)$,
$\Delta(k=1)$
\end{tabular}
\end{table}
Here $\mathbb B^n:=\{z\in \mathbb C^n: |z|<1\}$ denotes the $n$-dimensional unit ball.

\section{Totally geodesic isometric discs in polydiscs}
For a given $C^\infty$ map $f\colon M\to N$ between Riemannian manifolds $M$ and $N$, 
the differential $df$ is a section of $\text{Hom}(TM, E)$ where $E:= f^* TN$ is the pull-back bundle of $TN$ with respect to $f$ over $M$. 
Since $\text{Hom}(TM, E)$ is canonically identified with $E\otimes T^*M$, we may consider $df$ as an $E$-valued $1$-form on $M$.
Let $D^N$ denote the Levi-Civita connection on $N$ and $D:=f^* D^N$ be a naturally induced metric connection given on $E$. 
A map $f\colon M\to N$ is said to be totally geodesic if and only if for any vector fields $X$ on $M$, $D_X df=0$. Moreover, $f$ is totally geodesic if and only if $f$ maps any geodesic on $M$ onto a geodesic on $N$.

We say that a map $f\colon\Delta\to \Delta$ preserves geodesics if for any geodesic $\gamma$ in $\Delta$, $f\circ\gamma$ is a geodesic with possibly different speed.
Note that if $f$ preserves geodesics, then it is a local diffeomorphism.

\bl\label{unit disk-2}
Let $f\colon\Delta\to \Delta$ be a $C^2$-smooth map that preserves geodesics. Suppose that $f$ extends $C^1$ up to an open set $U\subset \partial\Delta$ and  $f_*(\mathbf{n}_p)$ is nowhere vanishing on $U$, where $\mathbf{n}_p$ is the outward unit normal vector at $p$. Then either $f$ or $\bar f$ is a holomorphic automorphism of $\Delta$.
\el 
\bpf
Choose a point in $\Delta$ and a point in $U$, say $0$ and $1$. After composing an automorphism of $\Delta$, we may assume that $f(0)=0$ and $f(1)=1$. Let $\gamma\colon\mathbb R\to \Delta$ be a unit speed geodesic such that $\gamma(0)=0$ and $\lim_{t\to \infty}\gamma(t)=1$. Then we obtain
$$\gamma(t)=\tanh(t).$$
Let $\tilde\gamma:=f\circ \gamma$.
Then we can write
$$\tilde\gamma(t)=\tanh(at),$$
for some $a>0$.
Suppose $a\neq 1$. Let 
$$\zeta=\tanh(  t),\quad \tanh(at)=h(\zeta).$$ 
Then the absolute value of 
$\lim_{\zeta\to 1}h'(\zeta)$ is either $0$ or $\infty$.
Since
$f(\zeta)=h(\zeta)$, for any $\zeta\in (-1,\,1)$
the absolute value of
$\lim_{\zeta\to 1} \partial_x f(\zeta)=f_*(\mathbf{ n}_1)$ is either $0$ or $\infty$,
contradicting the assumption on the radial derivative of $f$. Hence we obtain $a=1$ and as a result $f$ sends any unit speed geodesic emitted from $0$ toward a point in $U$ to a unit speed geodesic. By continuity $f$ sends any unit speed geodesic to a unit speed geodesic. This implies that $f$ is an isometry and hence it is an automorphism or a conjugate automorphism. 
\epf

\bp\label{isom-disc}
Let $f\colon (\Delta , \lambda g_{\Delta}) \to (\Delta^p, g_{\Delta^p})$ be a totally geodesic isometric embedding which extends $C^1$-smoothly over the boundary for $\lambda>0$. Then $\lambda$ is a positive integer and $f$ is of the form $f=(f_1, \dots f_\lambda, 0,\ldots, 0)$
where $f_j$ is of the form $\zeta\mapsto \zeta$ or $\zeta\mapsto \bar \zeta$ up to automorphisms. 
\ep

\bpf
Write $f=(f_1,\ldots, f_p)$. 
We may assume $f(0)=0$ and $f'(0) = (a_1,\ldots, a_k, 0,\ldots, 0)$ with nonzero $a_j$. Notice that $f_j\equiv 0$ for all $j=k+1,\ldots, p$.
First we will show the normal derivative $\frac{\partial f_j}{\partial {\mathbf n}}$ does not vanish identically on $\partial \Delta$, where ${\mathbf n}$ is the unit vector normal to $\partial \Delta$ for any nonconstant $f_j$.
Suppose that $\frac{ \partial f_j}{\partial {\mathbf n}}(p)=0$ for any $p\in \partial \Delta$.
Since each $f_j$ is totally geodesic, $f_j$ is a harmonic map with respect to metrics $\lambda g_\Delta$ and $g_\Delta$ for any $j$ (for detail see Example 6.3 in \cite{Wu88}). Let $u:= |f_j|^2-1  \colon (\Delta,\lambda g_\Delta)\to (-1,1)\subset (\Delta, g_\Delta)$ which is a subharmonic function, 
i.e. $\tau (u)\geq 0$ where $\tau (u)$ denotes the tension field of $u$.
Note that since $\frac{\partial f_j}{\partial\nu}\equiv 0$ on $\partial\Delta$, 
we have $u \equiv \frac{\partial u}{\partial \nu}\equiv 0$ on $\partial\Delta$. 
Since 
$$
0\leq \int_\Delta \tau(u) dV = \int_{\partial\Delta} \frac{\partial u}{\partial \nu} d\sigma =0
$$
by Green's formula, $u$ is a harmonic map.
Since $u \equiv \frac{\partial u}{\partial \nu}\equiv 0$ on $\partial\Delta$, $u$ is constant on $\Delta$ which is a contradiction.

Let $\gamma\colon \mathbb R\rightarrow \Delta$ be a complete geodesic given by $\gamma(t) = \tanh t$. Composing with an automorphism we may assume $\frac{\partial f_j}{\partial{\mathbf n}}(1)\neq 0$ for any $j=1,\ldots, k$ and 
$$f\circ \gamma (t) = (\tanh a_1 t, \tanh a_2 t,\ldots, \tanh a_k t, 0,\cdots,0)$$ for some nonzero $a_j\in \mathbb R$.
By the chain rule, we have 
\begin{equation}
\begin{aligned}
\frac{d}{dt} (f\circ \gamma)(t) 
=\frac{\partial \gamma}{\partial t} \left\{ \left( \frac{ \partial f_j}{\partial z} + \frac{ \partial f_j}{\partial \bar z} \right)\frac{\partial}{\partial z_j} 
+ \left( \frac{\partial \bar f_j}{\partial z} + \frac{\partial \bar f_j}{\partial \bar z}\right) \frac{\partial}{\partial \bar z_j}\right\}.
\end{aligned}
\end{equation}
Since $f$ is an isometric embedding, we have 
\begin{equation}\label{eq_iso}
\bigg| \frac{d\gamma}{dt}\bigg|^2 \sum_{j=1}^k \frac{
\big| \frac{ \partial f_j}{\partial z}  + \frac{ \partial f_j}{\partial \bar z} \big|^2 }{(1-\tanh^2 a_jt)^2}
=\sum_{j=1}^k \frac{ \cosh^2 a_jt
 }{\cosh^2 t}\bigg| \frac{ \partial f_j}{\partial z}  + \frac{ \partial f_j}{\partial \bar z} \bigg|^2
= \sum_{j=1}^k |a_j|^2 = \lambda.
\end{equation}
Since $f_1\colon\Delta\to\Delta$ preserves geodesics, by Lemma \ref{unit disk-2} we have $a_1=1$ and $f_1\in\Aut(\Delta)$ or $f_1\in\overline{\Aut(\Delta)}$. 
Moreover for any $\theta\in (-\varepsilon, \varepsilon)$ with sufficiently small $\varepsilon$, $f_1\circ (e^{i\theta}\tanh t)$ is also a unit speed geodesic and hence for any $\theta\in [-\pi, \pi]$, $f_1\circ (e^{i\theta}\tanh t)$ is a unit speed geodesic.
As a result, $f_1$ is an automorphism of $\Delta$.
By \eqref{eq_iso} we have 
\begin{equation}\nonumber 
\bigg| \frac{d\gamma}{dt}\bigg|^2 \sum_{j=2}^k \frac{
\big| \frac{ \partial f_j}{\partial z}  + \frac{ \partial f_j}{\partial \bar z} 
\big|^2  }{(1-\tanh^2 a_jt)^2}
= \lambda-1.
\end{equation}
By keeping to apply the above argument, we obtain that $k= \lambda$ and the lemma.
\epf
By the similar way, we obtain the following.
\bc
Let $\mathbb B^n = \{z\in \mathbb C^n : |z|<1\}$ be the $n$-dimensional unit ball. Any totally geodesic isometric embedding $f\colon \mathbb B^n\to\mathbb B^m$ is of the form $f(z_1,\ldots, z_n) = (z_1,\ldots, z_n, 0\ldots,0)$ or $f(z_1,\ldots, z_n) = (\bar z_1,\ldots, \bar z_n, 0\ldots,0)$ up to automorphisms.
\ec

%Let $p, q\in \Omega$. We call a holomorphic map $\varphi:\Delta\to \Omega$ {\em complex geodesic} joining $p$ and $q$ if $\varphi$ is a totally geodesic isometric embedding and $\varphi(\Delta)$ contains $p$ and $q$.

\section{Proof of Theorems}

\subsection{Proof of Theorem \ref{main thm}}

We will use induction on the rank of $\Omega$.
If rank$(\Omega)=1$, it is well-known (cf. \cite{GK82},\cite{ GS13},\cite{ A17}).
Assume that the theorem holds for bounded symmetric domains of rank $< r$ and assume that rank$(\Omega)=r$. 
%Assume that 
%the rank of $\dfrac{\partial f}{\partial x}(0)$ is the maximum among the rank of tangent vectors of $f(\Delta)$, where $x=Re(z)$ and 
%$$f(0)=F(0,0)=0\in \Omega.$$

\bl\label{main lem}
Let $f\colon\Delta\to \Omega\subset\mathbb C^n$ be a $C^2$ smooth totally geodesic isometric embedding that satisfies the condition in Theorem~\ref{main thm}. Then there exists a minimal disc $\Delta\subset\Omega$ passing through $f(0)$ such that
$$ f(\Delta)\subset \Delta\times\Delta^\perp$$
and each component of $f$ with respect to the above decomposition is a totally geodesic isometric embedding.
%where $(\Delta)^\perp$ is the characteristic subdomain of $\Omega$ passing through $f(0)$ orthogonal to $\Delta$.
\el

\bpf

Let $U\subset \partial\Delta$ be an open set where $f$ extends as a $C^1$ map. Fix a point $p\in U$. We may assume that $f_*(\mathbf n_p)\neq 0$.
Let $\ell$ be the rank of $f_*(\mathbf n_p)$. Choose the unique $\ell$-dimensional polysphere $P_p$ passing through $f(p)$ such that $f_*(\mathbf n_p)\in T_{f(p)}P_p$. 
For a point $x\in \Delta$, choose the unit speed geodesic $\gamma_x$ such that 
$\gamma_x(0)=x$ and $\lim_{t\to \infty}\gamma_x(t)=p$. Then 
we obtain
$$\lim_{t\to \infty}f_*(\mathbf n_{\gamma_x(t)})=f_*(\mathbf n_p),$$
where 
$$\mathbf n_z =\frac{1}{|z|}\left(x\frac{\partial}{\partial x}+y\frac{\partial }{\partial y}\right)$$
for $z=x+\sqrt{-1}y\neq 0$.
Since $f\circ\gamma_x$ is a geodesic in $\Omega$, we obtain
$f(\gamma_x)\subset P_p\times (P_p)^\perp$.
Since $x$ is arbitrary, we obtain $f(\Delta)\subset P_p\times (P_p)^\perp$.

Let $\Omega_1:=P_p\cap\Omega$ and $\Omega_2:=(P_p)^\perp\cap \Omega$. Then $f(\Delta)\subset \Omega_1\times\Omega_2$. 
Since $f$ is a totally geodesic isometric embedding, the $\Omega_1$ component of $f$, say $f_1$, preserves geodesics. Write $f_1=(g_1,\ldots,g_{\ell})\colon\Delta\to \Omega_1=(\Delta)^{\ell}$. Since $P_p$ is the minimal polysphere whose tangent space at $f(p)$ contains $f_*(\mathbf n_p)$, the radial derivative of each $g_i$ at $p$ is non-vanishing. Hence by Lemma~\ref{unit disk-2}, 
each $g_i$ is either holomorphic or anti holomorphic isometric map with respect to the Poincar\'{e} metric on the disk. Choose one of $g_i$'s. Then the conclusion follows.
\epf

\subsection{Proof of Theorem \ref{main cor}}
Assume that $f$ is neither holomorphic nor anti holomorphic.
%Suppose $f$ extends $C^1$ up to an open set $U$ of the Shilov boundary of $\Omega$ and the radial derivative of $f$ is nowhere vanishing on $U$.
%First we will show that $f:\Omega\to \Omega'$ is pluriharmonic.
%It is enough to show that $\partial\bar\partial f=0$.
Let $U$ be an open set in the Shilov boundary of $\Omega$ to which $f$ extends as a $C^1$ map.
Fix a point, say $0$, in $\Omega$. For a point $q\in U$, choose a totally geodesic maximal polysphere $P_q$ that contains $0$ and $q$. 
Let 
$$\Delta_q=\{\zeta q:\zeta\in \Delta\}\subset P_q\cap \Omega=:\Delta^r_q.$$ 
Then $\Delta_q$ is a totally geodesic disc with the largest holomorphic sectional curvature. Since $U\cap \partial\Delta_q$ is open in $\partial\Delta_q$, by Theorem~\ref{main thm}, $f$ is harmonic on $\Delta_q$. Since $U$ is open in the Shilov boundary of $\Omega$, the span of $T_0\Delta_q,~q\in U$ becomes $T_0\Omega$. Therefore we obtain $\partial\bar\partial f(0)=0$. Since $0$ is arbitrary, we obtain $\partial\bar\partial f\equiv 0, $ i.e. $f$ is pluriharmonic on $\Omega$. 
Therefore we can write 
$$f(z)=F(z)+G(\bar z)$$
for some nonconstant holomorphic maps $F$ and $G$.
Consider a holomorphic map $H\colon\Omega\times\Omega\to \mathbb C^{\dim \Omega'}$ defined by
$$ H(z, w)=F(z)+G(w), \quad (z, w)\in \Omega\times\Omega.$$
Then
$$H(z,\bar z)=f(z),\quad z\in \Omega.$$
We will show that 
$$H(\Omega\times\Omega)\subset \Omega'.$$
%and for any $w\in \Omega$,
%$H(\cdot,w):\Omega\to\Omega'$ and $H(w, \cdot):\Omega\to \Omega'$ are totally geodesic isometric embedding.

Let $0\in \Omega$, $q\in U$ and $\Delta^r_q$ be as above. For a point $x\in \Delta^r_q$, choose a complete geodesic $\gamma=(\gamma_1,\ldots,\gamma_r)\colon\mathbb R\to \Delta^r_q$ passing through $x$ such that each $\gamma_i$ is a unit speed geodesic in $\Delta$ and $\lim_{t\to\infty}\gamma(t)=q\in (\partial\Delta)^r_q$. Then there exists a totally geodesic disc $\Delta_\gamma\subset \Delta^r_q$ that contains the image of $\gamma$. Since $f$ is totally geodesic, by Theorem~\ref{main thm}, Lemma~\ref{isom-disc} and the argument in the proof of Lemma~\ref{main lem}  we can show that there exists a totally geodesic polysphere $P'_q$
passing through $f(q)$ and depending only on $f_*(\mathbf n_q)$ such that $f(\Delta_\gamma)\subset P'_q\cap \Omega'$. Since $x$ is arbitrary, we obtain $f(\Delta^r_q)\subset P'_q\cap\Omega'$.
Write $f=(f_1,\ldots,f_s)\colon\Delta^r_q\to P'_q\cap\Omega'=\Delta^s.$ By continuity of $f_*(\mathbf n_q)$, after shrinking $U$ if necessary, we may assume that $s$ is constant on $U$.

Consider a continuous family of totally geodesic discs 
$$\Delta_\alpha:=\{\zeta\alpha: \zeta\in \Delta\}\subset \Delta^r_q,\quad \alpha\in (\partial\Delta)^{r}_q\cap U.$$ 
Since $\Delta_\alpha$ is totally geodesic, by Theorem~\ref{main thm}, each $f_i$ is either holomorphic or anti holomorphic on $\Delta_\alpha$. Since we assumed that $f$ is neither holomorphic nor anti holomorphic, by continuity of derivatives of $f$, we may assume that 
$f_1,\ldots,f_{k}$ are holomorphic and $f_{k+1},\ldots,f_s$ are anti-holomorphic in $\Delta_\alpha$ for all $\alpha\in (\partial\Delta)^r_q\cap U$ for some integer $1\leq k<s$. Moreover, by continuity of derivatives of $f$, we may assume that $k$ is independent of $q\in U$. 
Since $U\cap (\partial\Delta)^r_q$ is open in the Shilov boundary of $\Delta^r_q$, the vector spaces $T_0\Delta_\alpha,~\alpha\in U\cap (\partial\Delta)_q^r$ span $T_0\Delta^r_q$ and therefore
we obtain that on $T_0\Delta^r_q,$
$$ \bar\partial f_i(0)=0,\quad i=1,\ldots,k$$
and
$$ \partial f_i(0)=0,\quad i=k+1,\ldots,s.$$
Similarly, we can show that
$$ \bar\partial f_i(x)=0,\quad i=1,\ldots,k$$
and
$$ \partial f_i(x)=0,\quad i=k+1,\ldots,s$$
for all $x\in \Delta^r_q$, i.e.
$$ f_i(z)=F_i(z),\quad i=1,\ldots,k$$
and
$$ f_i(z)=G_i(\bar z),\quad i=k+1,\ldots,s.$$
Therefore for $q\in U$, we obtain 
$$H(\Delta^r_q\times\Delta^r_q)=(f_1,\ldots,f_k)(\Delta^r_q)\times(f_{k+1},\ldots,f_s)(\Delta^r_q)\subset \Omega'$$
and the second fundamental forms of $F$, $G$ are trivial on $\Delta_q^r$ for $q\in U$. 
Furthermore, since each $F_i$, $G_i$ is an isometry on general totally geodesic disc with largest holomorphic sectional curvature, say $\Delta_d$, in $\Delta_q^r$, by Lemma~\ref{isom-disc}
the map $F_i|_{\Delta_d}\colon \Delta_d\to \Delta$ is an automorphism. 
Note that the induced metric on $\Delta_d$ from $\Omega$ is equivalent to $rg_\Delta$.
Hence 
$$
F|_{\Delta_d}^*(g_{\Delta^s}) 
= \sum_{i=1}^k F_i^* (g_{\Delta})
= kg_\Delta = \frac{k}{r} g_{\Delta_d}.
$$
Since such $\Delta_d$ are generic, by continuity we obtain
$$F^*(g_{\Delta^s})=\frac{k}{r}~g_{\Delta_q^r}
$$
By the similar way we obtain
$$
G^*(g_{\Delta^s})=\frac{s-k}{r}~g_{\Delta_q^r}.$$
This implies that 
$F$ and $G$ are totally geodesic isometric embedding on $\Delta^r_q$ with respect to the Bergman metric. Since $U$ is open in the Shilov boundary, the second fundamental forms of $F$ and $G$ are trivial on $\Omega$
and therefore $F$ and $G$ are totally geodesic isometric embedding with respect to the Fubini-Study metric of the compact duals. 
In particular, $F+G$ is $C^1$ up to the Shilov boundary of $\Omega\times\Omega$ and the radial derivative of $F+G$ is nontrivial there.
Since $f(z)=F(z)+G(\bar z)$, $f$ is $C^1$ up to the Shilov boundary of $\Omega$ and the radial derivative of $f$ is nonvanishing on a dense open set of the Shilov boundary, implying that
for all general $q\in S(\Omega)$, we obtain
$$H(\Delta^r_q\times\Delta^r_q)=(f_1,\ldots,f_k)(\Delta^r_q)\times(f_{k+1},\ldots,f_s)(\Delta^r_q)\subset \Omega'.$$
Fix $w=0$. Since $\Omega\subset \overline{~\bigcup_q \Delta^r_q~}$, where the union is taken over all general point $q\in S(\Omega)$,  by holomorphicity of $F$, we obtain
$$H(\Omega, 0)=F(\Omega)+G(0)\subset \bigcup_{q\in S(\Omega)}(f_1,\ldots,f_k)(\Delta^r_q)\times(f_{k+1},\ldots,f_s)(\Delta^r_q)\subset  \Omega'.$$
Since $0$ is an arbitrary point, we obtain
$$H(\Omega, w)=F(\Omega)+G(w)\subset   \Omega'$$
and by taking union over all $w\in \Omega$, we obtain
$$H(\Omega, \Omega)\subset  \Omega'.$$

Let $Y_1$ and $Y_2$ be the minimal totally geodesic subdomains of $\Omega'$ that contain $F(\Omega)+G(0)$ and $F(0)+G(\Omega)$, respectively.
We claim that $Y_1$ and $Y_2$ are totally geodesic in $\Omega'$ such that $Y_1\times Y_2$ can be embedded totally geodesically in $\Omega'$ and
$f=(f_1, f_2)\colon\Omega\to Y_1\times Y_2$, where $f_1(z)=F(z)$ and $f_2(z)=G(\bar z)$.
Choose a general totally geodesic disc $\Delta\subset \Omega$. Then as above, we can show that there exist polydiscs $P_1$, $P_2$ that are orthogonal to each other such that 
$$
F(\Delta)+G(0)\subset P_1
\quad\text{ and } \quad
F(0)+G(\Delta)\subset P_2.
$$ 
Furthermore, since $G$ is an isometry on $\Delta$, for $z_0\in \partial \Delta$, we obtain 
$$
F(0)+G(\bar z_0)\subset \partial P_2\subset\partial \Omega'.
$$ 
Since $\Omega'$ is pseudoconvex and $F$ is holomorphic, $F(\Delta)+G(\bar z_0)$ is contained in a boundary component $\Omega_{\bar z_0}'$
and by the same reason, $F(\Omega)+G(\bar z_0)$ is contained in the same boundary component $\Omega_{\bar z_0}'$. 
In particular, $F(\Omega)$ is orthogonal to $G(\Delta)$. Since $\Delta$ is arbitrary, $Y_1\times Y_2$ can be embedded totally geodesically into $\Omega'$.
Since $F$ is holomorphic and $G$ is anti-holomorphic, the conclusion follows.

\end{document}